\documentclass[preprint,12pt]{elsarticle}

\usepackage{amssymb}
\usepackage{amsmath}
\allowdisplaybreaks

\begin{document}

\begin{frontmatter}

\title{Formul\ae\ for the Number of\\ Partitions of $n$ into at most $m$ parts \\(Using the Quasi-Polynomial Ansatz)}

 \author[AVS]{Andrew V. Sills\fnref{fn1}}
 \ead{ASills@GeorgiaSouthern.edu}
 \ead[url]{http://math.georgiasouthern.edu/~asills}
 
 \author[DZ]{Doron Zeilberger\fnref{fn2}}
\ead{zeilberg@math.rutgers.edu}
\ead[url]{http://www.math.rutgers.edu/~zeilberg}

\fntext[fn1]{A. V. S. thanks DIMACS for hospitality during his July 2011 stay, in which this
research with D. Z. was initiated.}
\fntext[fn2]{D. Z. is partially supported by a grant from the NSF.}
 
  \address[AVS]{Department of Mathematical Sciences, 
Georgia Southern University, Statesboro, GA, 30458-8093, USA}

\address[DZ]{Department of Mathematics, Hill Center- Busch Campus, Rutgers University,
110 Frelinghuysen Rd, Piscataway, NJ 08854-8019, USA}

\begin{abstract}
The purpose of this short article is to announce, and briefly describe, a Maple package,
{\tt PARTITIONS}, that (inter alia) {\it completely automatically}  discovers, and then proves,
explicit expressions (as sums of quasi-polynomials) for $p_m(n)$ for any desired $m$.
We do this to demonstrate the power of ``rigorous guessing'' as facilitated by the
quasi-polynomial ansatz.
\end{abstract}


\end{frontmatter}

\section{Introduction}
Recall that a {\it partition} of a non-negative integer $n$ is a non-increasing sequence
of positive integers $\lambda_1 \dots \lambda_m$ that sum to $n$. For example
the integer $5$ has the following seven partitions: $\{5, 41, 32,311, 221, 2111, 11111\}$.
The {\bf bible} on partitions is George Andrews' {\it magnum opus} \cite{A1}.

We denote by $p_m(n)$ the number of partitions of $n$ into at most $m$ parts.
By a classic theorem~\cite[p. 8, Thm. 1.4]{A1}, $p_m(n)$ also equals the number of partitions of $n$ into parts that are at most $m$.
There is an 
extensive literature concerning formul\ae\ for $p_m(n)$,  including contributions
by Cayley, Sylvester, Glaisher, and Gupta.
 {For additional references and
historical notes,
see} George Andrews' 
fascinating article~\cite[\S3]{A2} {and Gupta's \textit{Tables}~\cite[pp. i--xxxix]{G}.
For an exhaustive history through 1920, see Dickson~\cite[Ch. 3]{D}. }

 More recently, George Andrews' student, Augustine O. Munagi, developed a beautiful theory
of so-called $q$-partial fractions~\cite{M}, {where the denominators
in the decomposition are always
expressions of the form $(1-q^r)^s$, rather than powers of cyclotomic polynomials as is
the case with the ordinary partial fraction decomposition.   Accordingly, formul\ae\  for
$p_m(n)$ derived from the $q$-partial fraction decomposition of the generating function
are most naturally expressed in terms of binomial coefficients.} 

It is well-known and easy to see that for any $m$, $p_m(n)$ is a {sum of}
 {\it quasi-polynomial{s of periods $1,2,3,\dots,m$}}.
A {\it quasi-polynomial of period $r$} is a function
{$f(n)$} on the integers
such that there exist $r$ polynomials $P_1(n), P_2(n), \dots, P_r(n)$ such that
$f(n)=P_i(n)$ if $n \equiv i\pmod{r}$. We represent such a quasi-polynomial as
a list $[P_1(n), \dots, P_r(n)]$.

Thus, e.g., we have, for $n\geq 0$, 
\begin{align} p_1(n) &=1,  \label{p1n} \\
p_2(n) &=\left[\frac n2 + \frac 34\right] + \left[-\frac 14, \frac 14 \right] ,\\
 %= \frac n2 + \frac 34 + \frac{(-1)^n}{4}
p_3(n) &= \left[ \frac{n^2}{12} + \frac{n}{2} + \frac{47}{72} \right] + \left[  -\frac 18, \frac 18 \right] +
 \left[ -\frac 19, -\frac 19, \frac 29 \right] \\
 p_4(n) &= \left[ \frac{n^3}{144} + \frac{5n^2}{48} + \frac{15n}{32} + \frac{175}{288} \right] +
 \left[ -\frac{n+5}{32},  \frac{n+5}{32} \right]   + \left[ 0, -\frac 19, \frac 19 \right] \notag
 \\ & \qquad \qquad+ \left[ 0, -\frac 18, 0, \frac 18\right]\\
 p_5(n) &= \left[ \frac{n^4}{2880} + \frac{n^3}{96} + \frac{31}{288}n^2 + \frac{85}{192} n + \frac{50651}{86400} \right]
   + \left[ -\frac{n}{64} - \frac{15}{128}, \frac{n}{64} + \frac{15}{128} \right] \notag \\ & \qquad\qquad
   + \left[ -\frac{1}{27}, -\frac{1}{27}, \frac{2}{27} \right] + \left[ \frac{1}{16}, -\frac{1}{16}, -\frac{1}{16}, \frac{1}{16}\right]
    \notag \\ & \qquad\qquad
   + \left[ -\frac{1}{25}, -\frac{1}{25}, -\frac{1}{25}, -\frac{1}{25}, \frac{4}{25} \right]. \label{p5n}
\end{align}

Eqs.~\eqref{p1n}--\eqref{p5n} were given in 1856 by Cayley~\cite[p. 132]{C} in a somewhat different form.  
In 1909, Glaisher~\cite{G1} presented formul\ae\ for $p_m(n)$ for $m=1,2,\dots,10$.  
In 1958, Gupta~\cite{G} extended Glaisher's results to the cases $m=11,12$.
In his 2005 Ph.D. thesis~\cite{Mt}, Munagi gave formul\ae\ for the cases $m=1,2,\dots,15$.
Munagi's formul\ae\ were derived with the aid of a Maple package he developed, and are of a 
somewhat different character than earlier contributions, as they follow from his
theory of $q$-partial fractions~\cite{M}.

\section{The \texttt{PARTITIONS} Maple package}
{\subsection{Overview}}
The purpose of this short article is to announce and briefly describe a Maple package,
{\tt PARTITIONS}, that {\it completely automatically}  discovers and proves
explicit expressions (as sums of quasi-polynomials) for $p_m(n)$ for any desired $m$. So far we only bothered
to derive the formul\ae\  for $1 \leq m \leq 70$, but one can easily go far beyond.

{\bf Not only that}, we can, more generally, derive (and prove!), completely automatically, 
expressions, as sums of quasi-polynomials,
for the number of ways of making change for $n$ cents in a country whose coins have
denominations of any given set of positive integers.

{\bf Not only that}, we can derive (and prove!), completely automatically, 
expressions (as sums of quasi-polynomials)
for $D_k(n)$, the number of partitions of $n$ whose {\it Durfee square}
has size $k$, for any desired, (numeric) positive integer $k$. 
(Recall that the size of the Durfee square of a partition $\lambda_1 \dots ...\lambda_m$
is the  largest $k$ such that $\lambda_k \geq k$.)

{\bf Not only that}, we (or rather our computers (and of course yours, if it has Maple and is loaded with our package))
can derive {\bf asymptotic expressions}, to {\it any desired order}, for both $p_m(n)$ and $D_k(n)$.
As far we we know the formula for $D_k(n)$ is brand-new, and the previous attempts for 
the asymptotic formula for $p_m(n)$ by humans G.J. Rieger~\cite{R} and E.M. Wright~\cite{W} (of Hardy-and-Wright fame)
only went as far as $O(n^{-2})$ and $O(n^{-4})$ respectively. We go all the way to $O(n^{-100})$!
(and of course can easily go far beyond).

{\bf Not only that}, 
we implement George Andrews' ingenious way \cite[sec. 3]{A2} to convert any quasi-polynomial to
a polynomial expression where one is also allowed to use the integer-part function
{$\lfloor n \rfloor$.}  This enabled our computers to find Andrews-style expressions for
$p_m(n)$ for $1 \leq m \leq 70$.

All these feats (and more!) are achieved by  the {\bf Maple package} {\tt PARTITIONS}. 

{\subsection{Using the PARTITIONS package}}
In order to use {PARTITIONS},
you must have {\tt Maple}$^{TM}$ installed on your computer. Then download the file:

{\tt http://www.math.rutgers.edu/{\char126\relax}zeilberg/tokhniot/PARTITIONS}
and save it as {\tt PARTITIONS}. Then launch Maple, and {at the
prompt,} enter:

{\tt read PARTITIONS}:

and follow the on-line instructions. Let's just highlight the most important
procedures.

{\tt AS100(m,n)}: 
shows the {\it pre-computed} first $100$ terms of the asymptotic
expression, in $n$, of $p_m(n)$ for {\it symbolic} $m$.

{\tt ASD80(k,n)}:
shows the {\it pre-computed} first $80$ terms of the asymptotic
expression, in $n$, of $D_k(n)$ for {\it symbolic} $k$.

{\tt BuildDBpmn(n,M)}: inputs a symbol $n$ and a {positive} integer $M$
and outputs a list of size $M$ whose $i$-th entry is an expression for $p_i(n)$ as a sum of $i$ quasi-polynomials

{\tt DiscoverAS(m,n,k)}: discovers the asymptotic expansion to order $k$
of $p_m(n)$ (the number of partitions of $n$ into at most $m$ parts) for large $n$ and fixed, but {\it symbolic}, $m$.

{\tt DiscoverDAS(k,n,r)}: discovers the asymptotic expansion to order $r$
of $D_k(n)$ (the number of partitions of $n$ whose Durfee square has size $k$)
for large $n$ and fixed, but {\it symbolic} $k$.

{\tt Durfee(k,n)}: discovers (rigorously!) the quasi-polynomial expression, in $n$, for $D_k(n)$,
for any desired positive integer $k$. It is extremely fast for small $k$, but of course gets slower as $k$ gets larger.

{\tt DurfeePC(k,n)}: does the same thing (much faster, of course!) using the pre-computed expressions of {\tt Durfee(k,n);}
for $k \leq 40$.

{\tt EvalQPS(L,n,n0)}: evaluates the sum of the quasi-polynomials in the variable $n$ 
given in the list $L$  at $n=n_0$.

{\tt HRR(n,T)}: evaluates in floating point the sum of the first $T$ terms of the Hardy-Ramanujan-Rademacher formula
for $p(n)$, the number of unrestricted partitions of $n$: \begin{equation*} \label{hrr}
     p(n) = \frac{1}{\pi\sqrt{2}} \sum_{k\geq 1} \sqrt{k} \ \underset{\gcd(h,k)=1}{ \sum_{0\leq h<k}} e^{\pi i \left( s(h,k)-2 n h/k  \right) }  \frac{d}{dn} 
   \left( \frac{\sinh\left( \frac{\pi}{k} \sqrt{ \frac23\left( n - \frac{1}{24}\right)}\right)}{\sqrt{n - \frac{1}{24}}} \right), 
\end{equation*}
where 
 $s(h,k) = \sum_{j=1}^{k-1} \left( \frac jk - \lfloor \frac jk \rfloor - \frac 12\right) \left( \frac{hj}k - \lfloor \frac{hj}k \rfloor - \frac 12\right)  $is the Dedekind sum.

Please be warned that for larger $n$ you need to increase {\tt Digits}. In order to get reliable results
you may want to use procedure {\tt HRRr(n,T,k)}.

{\tt pmn(m,n)}: discovers (rigorously!) the quasi-polynomial expression, in $n$, for $p_m(n)$,
for any desired positive integer $m$. It is extremely fast for small $m$, but of course gets slower as $m$ gets larger.

{\tt pmnPC(m,n)}: does the same thing (much faster, of course!) using the pre-computed expressions of {\tt pmn(m,n);}
for $m \leq 70$.

{\tt pmnAndrews(m,n)}: discovers (rigorously!) the Andrews-style expression, in $n$, for $p_m(n)$
for any desired positive integer $m$. Instead of using quasi-polynomials explicitly (that some humans
find awkward), it uses the integer-part function {$\lfloor n \rfloor$}, denoted by {\tt trunc(n)} in Maple.

{\tt pn(n)}:  the number of partitions of $n$, $p(n)$, using Euler's recurrence. It is useful for checking,
since $p_n(n)=p(n)$.

{\tt pnSeq(N)}:  the list of the first $N$ values of $p(n)$. The output of {\tt pnSeq(50000):}
can be gotten from

{\tt http://www.math.rutgers.edu/{\char126\relax}zeilberg/tokhniot/oPARTITIONS9}
where this list of 50000 terms is called {\tt pnTable}.

{\tt pSn(S,n,K)}: the more general problem where the parts are drawn from the list $S$ of positive integers.
It outputs an explicit expression, as a sum of quasi-polynomials,
for $p_S(n)$, the number of integer partitions of $n$ whose parts are drawn from the finite list of 
positive integers $S$. $K$ is a guessing parameter, that should be made higher if the procedure returns {\tt FAIL}.

{\tt pmnNum(m,n0)}: like {\tt pmn(m,n);} but for both {\it numeric} $m$ and $n0$. The output is a number.
For $m \leq 70$ it is extremely fast, since it uses the {\it pre-computed} values of $p_m(n)$ gotten from
{\tt pmnPC(m,n);}. For example to get the number of integer partitions of a googol ($10^{100}$) into at most
$60$ parts, you would get, in $0.02$ seconds, the $5783$-digit integer, by simply typing 

{\tt pmnNum(60,10**100);} \quad .

One of us (DZ) posed this is a $100$-dollar challenge to the users of the very useful {\tt Mathoverflow} forum.
This was taken-up, successfully,
by user {\it joro}\cite{J}, whose computer did it correctly in about $2$ hours, using {\tt PARI}.
User {\tt joro} generously suggested that instead of sending him a check, DZ should donate it
in {\tt joro}'s honor, to a charity of DZ's choice, and the latter decided on the Wikipedia Foundation.

{\bf Sample input and output} can be gotten from the ``front'' of this article:
\hfill\break
{\tt http://www.math.rutgers.edu/{\char126\relax}zeilberg/mamarim/mamarimhtml/pmn.html} \quad .

\section{Methodology: Rigorous Guessing}

{The idea of deriving formul\ae\ for $p_m(n)$ and
$p_S(n)$ with the aid of the partial fraction decomposition of the generating function dates} back {at least} to Cayley {\cite{C}}.  
We ask Maple to convert the generating function
\[ {\sum_{n\geq 0} p_m(n) q^n = \frac{1}{(1-q)(1-q^2) \cdots (1-q^m)}}\]
or in the case of $p_S(n)${, where $S=\{ s_1, s_2, \dots,s_j\}$,}
 \[ {\sum_{n\geq 0} p_S(n)q^n = \frac{1}{(1-q^{s_1})(1-q^{s_2}) \cdots (1-q^{s_j})} }\]
into {partial fractions}.
Then for each piece, Maple finds the first few terms of the Maclaurin expansion, and then fits
the data with an appropriate quasi-polynomial using {\it undetermined} coefficients. 
The output is the list of these quasi-polynomials whose sum is the desired expression
for $p_m(n)$ or $p_S(n)$.
See the source-code for more details.

\textbf{Example.}
{Consider the case $m=4$.  We have Maple calculate that
\begin{multline} \sum_{n\geq 0} p_4(n) q^n = \frac{1}{(1-q)(1-q^2)(1-q^3)(1-q^4)}
\\= \frac{17/72}{1-q} + \frac{59/288}{(1-q)^2} + \frac{1/8}{(1-q)^3} + \frac{1/24}{(1-q)^4}
+ \frac{1/8}{1+q} + \frac{1/32}{(1+q)^2} + \frac{(1+q)/9}{1+q+q^2}. \label{pfd}
\end{multline}
At this point we could, as Cayley did, expand each term as a series in $q$, collect like terms,
and then the coefficient of $q^n$ will be a formula for $p_4(n)$, but why bother?
From Sylvester~\cite{S} and Glaisher~\cite{G2}, we know that 
   \[ p_4(n) = \sum_{j=1}^4 W_j(n), \] where
 each $W_j(n)$ is a quasi-polynomial $[P_{j1}(n), P_{j2}(n), \dots, P_{jj}(n)]$ of period $j$.
Further, $W_j(n)$ is of degree $\lfloor \frac{m-j}{j} \rfloor$, and arises from those terms of~\eqref{pfd} with denominator
a power of the $j$-th cyclotomic polynomial.  
Instead, let us allow Maple to guess the correct quasi-polynomials:
  We know \textit{a priori} that $W_1(n)$ is of the form $[ a_0 + a_1 n + a_2 n^2 + a_3 n^3]$
and let Maple calculate the (beginning of the) Maclaurin series for the terms of~\eqref{pfd} that
contribute to $W_1(n)$:
  \[  \frac{17/72}{1-q} + \frac{59/288}{(1-q)^2} + \frac{1/8}{(1-q)^3} + \frac{1/24}{(1-q)^4}
   = \frac{175}{288} + \frac{19}{16}q + \frac{581}{288} q^2 + \frac{113}{36}q^3 + O(q^4).\]
  Thus,
  \[ \left[ {\begin{array}{cccc} 1 & 0 & 0 & 0 \\
                                                  1^0 & 1^1 & 1^2 & 1^3 \\
                                                  2^0 & 2^1 & 2^2  & 2^3 \\
                                                   3^0 & 3^1 & 3^2 & 3^3  \end{array}}\right]
   \left[ \begin{array}{c} a_0 \\ a_1 \\ a_2 \\a_3 \end{array} \right] =
   \left[ \begin{array}{c} {175}/{288} \\ {19}/{16} \\ {581}/{288} \\ {113}/{36} \end{array}\right],                                               
  \] 
 which immediately implies that 
 \[ W_1(n) = \left[\frac{1}{144} n^3 + \frac{5}{48} n^2 + \frac{15}{32} n + \frac{175}{288}\right]. \]
  Similarly, for $W_2(n)$, which must be of the form
  \[ [a_1 + a_3 n, a_0 + a_2 n  ] , \] we find
  \[ \frac{1/8}{1+q} + \frac{1/32}{(1+q)^2} = \frac{5}{32} - \frac{3}{16}q + \frac{7}{32}q^2 - \frac{1}{4} q^3 + O(q^4),\]
  so that
  \[ \left[ \begin{array}{cc} 1 & 0 \\ 1 & 2 \end{array} \right] 
     \left[ \begin{array}{c} a_0 \\ a_2 \end{array}  \right] 
    =\left[ \begin{array}{c}  5/32 \\ 7/32  \end{array}  \right]
  \]
 and
  \[ \left[ \begin{array}{cc} 1 & 1 \\ 1 & 3 \end{array} \right] 
     \left[ \begin{array}{c} a_1 \\ a_3 \end{array}  \right] 
    =\left[ \begin{array}{c}  -3/16 \\ -1/4  \end{array}  \right],
  \] and thus
  \[ W_2(n) = \left[ -\frac{5}{32} - \frac{n}{32} ,  \frac{5}{32} + \frac{n}{32} \right] . \]
  Analogous reasoning yields $W_3(n)=\left[ 0, -\frac19,\frac19\right]$ and 
  $W_4(n)=\left[ 0,-\frac 18,0,\frac18 \right]$.
}

\section{Conclusion}
The present approach uses \textit{very na\"ive} guessing to discover, and \textit{prove} (rigorously!),
formulas (or as Cayley and Sylvester would say, \textit{formul\ae}) for the number of partitions of
the integer $n$ into at most parts $m$ parts for $m\leq 70$, and of course, one can easily go far 
beyond.  The core of the idea goes back to Arthur Cayley, and is familiar to any second-semester
calculus student: partial fractions!  But dear Arthur could only go so far, so his good buddy, 
James Joseph Sylvester, designated a sophisticated theory of ``waves''~\cite{S} that facilitated hand
calculations, which were later dutifully carried out by J. W. L. Glaisher in~\cite{G2}.
But, with modern computer algebra systems (Maple in our case), one can go much further just
using Cayley's original ideas.

\section*{Acknowledgment}
The authors thank Ken Ono for several helpful comments on an earlier version of this manuscript.

\bibliographystyle{elsarticle-num}

\end{document}